\input amssym.def

\catcode`\@=11
\font\tensmc=cmcsc10      
\def\smc{\tensmc}
\def\pagewidth#1{\hsize= #1 }
\def\pageheight#1{\vsize= #1 }
\def\hcorrection#1{\advance\hoffset by #1 }
\def\vcorrection#1{\advance\voffset by #1 }
\def\wlog#1{}
\newif\iftitle@
\outer\def\title{\title@true\vglue 24\p@ plus 12\p@ minus 12\p@
   \bgroup\let\\=\cr\tabskip\centering
   \halign to \hsize\bgroup\tenbf\hfill\ignorespaces##\unskip\hfill\cr}
\def\endtitle{\cr\egroup\egroup\vglue 18\p@ plus 12\p@ minus 6\p@}
\outer\def\author{\iftitle@\vglue -18\p@ plus -12\p@ minus -6\p@\fi\vglue
    12\p@ plus 6\p@ minus 3\p@\bgroup\let\\=\cr\tabskip\centering
    \halign to \hsize\bgroup\smc\hfill\ignorespaces##\unskip\hfill\cr}
\def\endauthor{\cr\egroup\egroup\vglue 18\p@ plus 12\p@ minus 6\p@}
\outer\def\heading{\bigbreak\bgroup\let\\=\cr\tabskip\centering
    \halign to \hsize\bgroup\smc\hfill\ignorespaces##\unskip\hfill\cr}
\def\endheading{\cr\egroup\egroup\nobreak\medskip}
\outer\def\subheading#1{\medbreak\noindent{\tenbf\ignorespaces
      #1\unskip.\enspace}\ignorespaces}
\outer\def\proclaim#1{\medbreak\noindent\smc\ignorespaces
    #1\unskip.\enspace\sl\ignorespaces}
\outer\def\endproclaim{\par\ifdim\lastskip<\medskipamount\removelastskip
  \penalty 55 \fi\medskip\rm}
\outer\def\demo#1{\par\ifdim\lastskip<\smallskipamount\removelastskip
    \smallskip\fi\noindent{\smc\ignorespaces#1\unskip:\enspace}\rm
      \ignorespaces}
\outer\def\enddemo{\par\smallskip}
\newcount\footmarkcount@
\footmarkcount@=1
\def\makefootnote@#1#2{\insert\footins{\interlinepenalty=100
  \splittopskip=\ht\strutbox \splitmaxdepth=\dp\strutbox 
  \floatingpenalty=\@MM
  \leftskip=\z@\rightskip=\z@\spaceskip=\z@\xspaceskip=\z@
  \noindent{#1}\footstrut\rm\ignorespaces #2\strut}}
\def\footnote{\let\@sf=\empty\ifhmode\edef\@sf{\spacefactor
   =\the\spacefactor}\/\fi\futurelet\next\footnote@}
\def\footnote@{\ifx"\next\let\next\footnote@@\else
    \let\next\footnote@@@\fi\next}
\def\footnote@@"#1"#2{#1\@sf\relax\makefootnote@{#1}{#2}}
\def\footnote@@@#1{$^{\number\footmarkcount@}$\makefootnote@
   {$^{\number\footmarkcount@}$}{#1}\global\advance\footmarkcount@ by 1 }

\hyphenation{man-u-script man-u-scripts ap-pen-dix ap-pen-di-ces}
\hyphenation{data-base data-bases}
\ifx\amstexloaded@\relax\catcode`\@=13 
  \endinput\else\let\amstexloaded@=\relax\fi
\newlinechar=`\^^J
\def\eat@#1{}
\def\Space@.{\futurelet\Space@\relax}
\Space@. %
\newhelp\athelp@
{Only certain combinations beginning with @ make sense to me.^^J
Perhaps you wanted \string\@\space for a printed @?^^J
I've ignored the character or group after @.}
\def\futureletnextat@{\futurelet\next\at@}
{\catcode`\@=\active
\lccode`\Z=`\@ \lowercase
{\gdef@{\expandafter\csname futureletnextatZ\endcsname}
\expandafter\gdef\csname atZ\endcsname
   {\ifcat\noexpand\next a\def\next{\csname atZZ\endcsname}\else
   \ifcat\noexpand\next0\def\next{\csname atZZ\endcsname}\else
    \def\next{\csname atZZZ\endcsname}\fi\fi\next}
\expandafter\gdef\csname atZZ\endcsname#1{\expandafter
   \ifx\csname #1Zat\endcsname\relax\def\next
     {\errhelp\expandafter=\csname athelpZ\endcsname
      \errmessage{Invalid use of \string@}}\else
       \def\next{\csname #1Zat\endcsname}\fi\next}
\expandafter\gdef\csname atZZZ\endcsname#1{\errhelp
    \expandafter=\csname athelpZ\endcsname
      \errmessage{Invalid use of \string@}}}}
\def\atdef@#1{\expandafter\def\csname #1@at\endcsname}
\newhelp\defahelp@{If you typed \string\define\space cs instead of
\string\define\string\cs\space^^J
I've substituted an inaccessible control sequence so that your^^J
definition will be completed without mixing me up too badly.^^J
If you typed \string\define{\string\cs} the inaccessible control sequence^^J
was defined to be \string\cs, and the rest of your^^J
definition appears as input.}
\newhelp\defbhelp@{I've ignored your definition, because it might^^J
conflict with other uses that are important to me.}
\def\define{\futurelet\next\define@}
\def\define@{\ifcat\noexpand\next\relax
  \def\next{\define@@}%
  \else\errhelp=\defahelp@
  \errmessage{\string\define\space must be followed by a control 
     sequence}\def\next{\def\garbage@}\fi\next}
\def\undefined@{}
\def\preloaded@{}    
\def\define@@#1{\ifx#1\relax\errhelp=\defbhelp@
   \errmessage{\string#1\space is already defined}\def\next{\def\garbage@}%
   \else\expandafter\ifx\csname\expandafter\eat@\string
         #1@\endcsname\undefined@\errhelp=\defbhelp@
   \errmessage{\string#1\space can't be defined}\def\next{\def\garbage@}%
   \else\expandafter\ifx\csname\expandafter\eat@\string#1\endcsname\relax
     \def\next{\def#1}\else\errhelp=\defbhelp@
     \errmessage{\string#1\space is already defined}\def\next{\def\garbage@}%
      \fi\fi\fi\next}
\def\famzero{\fam\z@}

\def\exp{\mathop{\famzero exp}\nolimits}

\def\inf{\mathop{\famzero inf}}

\def\lim{\mathop{\famzero lim}}
\def\liminf{\mathop{\famzero lim\,inf}}

\def\log{\mathop{\famzero log}\nolimits}
\def\max{\mathop{\famzero max}}

\def\textfont@#1#2{\def#1{\relax\ifmmode
    \errmessage{Use \string#1\space only in text}\else#2\fi}}
\textfont@\rm\tenrm
\textfont@\it\tenit
\textfont@\sl\tensl
\textfont@\bf\tenbf
\textfont@\smc\tensmc
\let\ic@=\/
\def\/{\unskip\ic@}
\def\textfonti{\the\textfont1 }
\def\t#1#2{{\edef\next{\the\font}\textfonti\accent"7F \next#1#2}}
\let\B=\=
\let\D=\.
\def~{\unskip\nobreak\ \ignorespaces}
{\catcode`\@=\active
\gdef\@{\char'100 }}
\atdef@-{\leavevmode\futurelet\next\athyph@}
\def\athyph@{\ifx\next-\let\next=\athyph@@
  \else\let\next=\athyph@@@\fi\next}
\def\athyph@@@{\hbox{-}}
\def\athyph@@#1{\futurelet\next\athyph@@@@}
\def\athyph@@@@{\if\next-\def\next##1{\hbox{---}}\else
    \def\next{\hbox{--}}\fi\next}
\def\.{.\spacefactor=\@m}
\atdef@.{\null.}
\atdef@,{\null,}
\atdef@;{\null;}
\atdef@:{\null:}
\atdef@?{\null?}
\atdef@!{\null!}   
\def\srdr@{\thinspace}                     
\def\drsr@{\kern.02778em}
\def\sldl@{\kern.02778em}
\def\dlsl@{\thinspace}
\atdef@"{\unskip\futurelet\next\atqq@}
\def\atqq@{\ifx\next\Space@\def\next. {\atqq@@}\else
         \def\next.{\atqq@@}\fi\next.}
\def\atqq@@{\futurelet\next\atqq@@@}
\def\atqq@@@{\ifx\next`\def\next`{\atqql@}\else\def\next'{\atqqr@}\fi\next}
\def\atqql@{\futurelet\next\atqql@@}
\def\atqql@@{\ifx\next`\def\next`{\sldl@``}\else\def\next{\dlsl@`}\fi\next}
\def\atqqr@{\futurelet\next\atqqr@@}
\def\atqqr@@{\ifx\next'\def\next'{\srdr@''}\else\def\next{\drsr@'}\fi\next}
\def\flushpar{\par\noindent}
\def\textfontii{\the\textfont2 }
\def\{{\relax\ifmmode\lbrace\else
    {\textfontii f}\spacefactor=\@m\fi}
\def\}{\relax\ifmmode\rbrace\else
    \let\@sf=\empty\ifhmode\edef\@sf{\spacefactor=\the\spacefactor}\fi
      {\textfontii g}\@sf\relax\fi}   
\def\nonhmodeerr@#1{\errmessage
     {\string#1\space allowed only within text}}
\def\linebreak{\relax\ifhmode\unskip\break\else
    \nonhmodeerr@\linebreak\fi}
\def\allowlinebreak{\relax
   \ifhmode\allowbreak\else\nonhmodeerr@\allowlinebreak\fi}
\newskip\saveskip@
\def\nolinebreak{\relax\ifhmode\saveskip@=\lastskip\unskip
  \nobreak\ifdim\saveskip@>\z@\hskip\saveskip@\fi
   \else\nonhmodeerr@\nolinebreak\fi}
\def\newline{\relax\ifhmode\null\hfil\break
    \else\nonhmodeerr@\newline\fi}
\def\nonmathaerr@#1{\errmessage
     {\string#1\space is not allowed in display math mode}}
\def\nonmathberr@#1{\errmessage{\string#1\space is allowed only in math mode}}
\def\mathbreak{\relax\ifmmode\ifinner\break\else
   \nonmathaerr@\mathbreak\fi\else\nonmathberr@\mathbreak\fi}
\def\nomathbreak{\relax\ifmmode\ifinner\nobreak\else
    \nonmathaerr@\nomathbreak\fi\else\nonmathberr@\nomathbreak\fi}
\def\allowmathbreak{\relax\ifmmode\ifinner\allowbreak\else
     \nonmathaerr@\allowmathbreak\fi\else\nonmathberr@\allowmathbreak\fi}
\def\pagebreak{\relax\ifmmode
   \ifinner\errmessage{\string\pagebreak\space
     not allowed in non-display math mode}\else\postdisplaypenalty-\@M\fi
   \else\ifvmode\penalty-\@M\else\edef\spacefactor@
       {\spacefactor=\the\spacefactor}\vadjust{\penalty-\@M}\spacefactor@
        \relax\fi\fi}
\def\nopagebreak{\relax\ifmmode
     \ifinner\errmessage{\string\nopagebreak\space
    not allowed in non-display math mode}\else\postdisplaypenalty\@M\fi
    \else\ifvmode\nobreak\else\edef\spacefactor@
        {\spacefactor=\the\spacefactor}\vadjust{\penalty\@M}\spacefactor@
         \relax\fi\fi}
\def\newpage{\relax\ifvmode\vfill\penalty-\@M\else\nonvmodeerr@\newpage\fi}
\def\nonvmodeerr@#1{\errmessage
    {\string#1\space is allowed only between paragraphs}}
\def\smallpagebreak{\relax\ifvmode\smallbreak
      \else\nonvmodeerr@\smallpagebreak\fi}
\def\medpagebreak{\relax\ifvmode\medbreak
       \else\nonvmodeerr@\medpagebreak\fi}
\def\bigpagebreak{\relax\ifvmode\bigbreak
      \else\nonvmodeerr@\bigpagebreak\fi}
\newdimen\captionwidth@
\captionwidth@=\hsize
\advance\captionwidth@ by -1.5in
\def\caption#1{}
\def\topspace#1{\gdef\thespace@{#1}\ifvmode\def\next
    {\futurelet\next\topspace@}\else\def\next{\nonvmodeerr@\topspace}\fi\next}
\def\topspace@{\ifx\next\Space@\def\next. {\futurelet\next\topspace@@}\else
     \def\next.{\futurelet\next\topspace@@}\fi\next.}
\def\topspace@@{\ifx\next\caption\let\next\topspace@@@\else
    \let\next\topspace@@@@\fi\next}
 \def\topspace@@@@{\topinsert\vbox to 
       \thespace@{}\endinsert}
\def\topspace@@@\caption#1{\topinsert\vbox to
    \thespace@{}\nobreak
      \smallskip
    \setbox\z@=\hbox{\noindent\ignorespaces#1\unskip}%
   \ifdim\wd\z@>\captionwidth@
   \centerline{\vbox{\hsize=\captionwidth@\noindent\ignorespaces#1\unskip}}%
   \else\centerline{\box\z@}\fi\endinsert}
\def\midspace#1{\gdef\thespace@{#1}\ifvmode\def\next
    {\futurelet\next\midspace@}\else\def\next{\nonvmodeerr@\midspace}\fi\next}
\def\midspace@{\ifx\next\Space@\def\next. {\futurelet\next\midspace@@}\else
     \def\next.{\futurelet\next\midspace@@}\fi\next.}
\def\midspace@@{\ifx\next\caption\let\next\midspace@@@\else
    \let\next\midspace@@@@\fi\next}
 \def\midspace@@@@{\midinsert\vbox to 
       \thespace@{}\endinsert}
\def\midspace@@@\caption#1{\midinsert\vbox to
    \thespace@{}\nobreak
      \smallskip
      \setbox\z@=\hbox{\noindent\ignorespaces#1\unskip}%
      \ifdim\wd\z@>\captionwidth@
    \centerline{\vbox{\hsize=\captionwidth@\noindent\ignorespaces#1\unskip}}%
    \else\centerline{\box\z@}\fi\endinsert}
\mathchardef\prime@="0230
\def\prime{{{}\prime@{}}}
\def\prim@s{\prime@\futurelet\next\pr@m@s}

\def\,{\relax\ifmmode\mskip\thinmuskip\else\thinspace\fi}
\def\!{\relax\ifmmode\mskip-\thinmuskip\else\negthinspace\fi}
\def\frac#1#2{{#1\over#2}}

\def\:{\nobreak\hskip.1111em{:}\hskip.3333em plus .0555em\relax}
\def\intic@{\mathchoice{\hskip5\p@}{\hskip4\p@}{\hskip4\p@}{\hskip4\p@}}
\def\negintic@
 {\mathchoice{\hskip-5\p@}{\hskip-4\p@}{\hskip-4\p@}{\hskip-4\p@}}
\def\intkern@{\mathchoice{\!\!\!}{\!\!}{\!\!}{\!\!}}
\def\intdots@{\mathchoice{\cdots}{{\cdotp}\mkern1.5mu
    {\cdotp}\mkern1.5mu{\cdotp}}{{\cdotp}\mkern1mu{\cdotp}\mkern1mu
      {\cdotp}}{{\cdotp}\mkern1mu{\cdotp}\mkern1mu{\cdotp}}}
\newcount\intno@             
\def\iint{\intno@=\tw@\futurelet\next\ints@} 
\def\iiint{\intno@=\thr@@\futurelet\next\ints@}
\def\iiiint{\intno@=4 \futurelet\next\ints@}
\def\idotsint{\intno@=\z@\futurelet\next\ints@}
\def\ints@{\findlimits@\ints@@}
\newif\iflimtoken@
\newif\iflimits@
\def\findlimits@{\limtoken@false\limits@false\ifx\next\limits
 \limtoken@true\limits@true\else\ifx\next\nolimits\limtoken@true\limits@false
    \fi\fi}
\def\multintlimits@{\intop\ifnum\intno@=\z@\intdots@
  \else\intkern@\fi
    \ifnum\intno@>\tw@\intop\intkern@\fi
     \ifnum\intno@>\thr@@\intop\intkern@\fi\intop}
\def\multint@{\int\ifnum\intno@=\z@\intdots@\else\intkern@\fi
   \ifnum\intno@>\tw@\int\intkern@\fi
    \ifnum\intno@>\thr@@\int\intkern@\fi\int}
\def\ints@@{\iflimtoken@\def\ints@@@{\iflimits@
   \negintic@\mathop{\intic@\multintlimits@}\limits\else
    \multint@\nolimits\fi\eat@}\else
     \def\ints@@@{\multint@\nolimits}\fi\ints@@@}
\def\Sb{_\bgroup\vspace@
        \baselineskip=\fontdimen10 \scriptfont\tw@
        \advance\baselineskip by \fontdimen12 \scriptfont\tw@
        \lineskip=\thr@@\fontdimen8 \scriptfont\thr@@
        \lineskiplimit=\thr@@\fontdimen8 \scriptfont\thr@@
        \Let@\vbox\bgroup\halign\bgroup \hfil$\scriptstyle
            {##}$\hfil\cr}
\def\endSb{\crcr\egroup\egroup\egroup}
\def\Sp{^\bgroup\vspace@
        \baselineskip=\fontdimen10 \scriptfont\tw@
        \advance\baselineskip by \fontdimen12 \scriptfont\tw@
        \lineskip=\thr@@\fontdimen8 \scriptfont\thr@@
        \lineskiplimit=\thr@@\fontdimen8 \scriptfont\thr@@
        \Let@\vbox\bgroup\halign\bgroup \hfil$\scriptstyle
            {##}$\hfil\cr}
\def\endSp{\crcr\egroup\egroup\egroup}
\def\Let@{\relax\iffalse{\fi\let\\=\cr\iffalse}\fi}
\def\vspace@{\def\vspace##1{\noalign{\vskip##1 }}}
\def\aligned{\,\vcenter\bgroup\vspace@\Let@\openup\jot\m@th\ialign
  \bgroup \strut\hfil$\displaystyle{##}$&$\displaystyle{{}##}$\hfil\crcr}
\def\endaligned{\crcr\egroup\egroup}
\def\matrix{\,\vcenter\bgroup\Let@\vspace@
    \normalbaselines
  \m@th\ialign\bgroup\hfil$##$\hfil&&\quad\hfil$##$\hfil\crcr
    \mathstrut\crcr\noalign{\kern-\baselineskip}}
\def\endmatrix{\crcr\mathstrut\crcr\noalign{\kern-\baselineskip}\egroup
                \egroup\,}
\newtoks\hashtoks@
\hashtoks@={#}
\def\format{\crcr\egroup\iffalse{\fi\ifnum`}=0 \fi\format@}
\def\format@#1\\{\def\preamble@{#1}%
  \def\c{\hfil$\the\hashtoks@$\hfil}%
  \def\r{\hfil$\the\hashtoks@$}%
  \def\l{$\the\hashtoks@$\hfil}%
  \setbox\z@=\hbox{\xdef\Preamble@{\preamble@}}\ifnum`{=0 \fi\iffalse}\fi
   \ialign\bgroup\span\Preamble@\crcr}

\def\cases{\left\{\,\vcenter\bgroup\vspace@
     \normalbaselines\openup\jot\m@th
       \Let@\ialign\bgroup$##$\hfil&\quad$##$\hfil\crcr
      \mathstrut\crcr\noalign{\kern-\baselineskip}}

\newif\iftagsleft@
\tagsleft@true
\def\TagsOnRight{\global\tagsleft@false}
\def\tag#1$${\iftagsleft@\leqno\else\eqno\fi
 \hbox{\def\pagebreak{\global\postdisplaypenalty-\@M}%
 \def\nopagebreak{\global\postdisplaypenalty\@M}\rm(#1\unskip)}%
  $$\postdisplaypenalty\z@\ignorespaces}
\interdisplaylinepenalty=\@M
\def\allowdisplaybreak@{\def\allowdisplaybreak{\noalign{\allowbreak}}}
\def\displaybreak@{\def\displaybreak{\noalign{\break}}}
\def\align#1\endalign{\def\tag{&}\vspace@\allowdisplaybreak@\displaybreak@
  \iftagsleft@\lalign@#1\endalign\else
   \ralign@#1\endalign\fi}
\def\ralign@#1\endalign{\displ@y\Let@\tabskip\centering\halign to\displaywidth
     {\hfil$\displaystyle{##}$\tabskip=\z@&$\displaystyle{{}##}$\hfil
       \tabskip=\centering&\llap{\hbox{(\rm##\unskip)}}\tabskip\z@\crcr
             #1\crcr}}
\def\lalign@
 #1\endalign{\displ@y\Let@\tabskip\centering\halign to \displaywidth
   {\hfil$\displaystyle{##}$\tabskip=\z@&$\displaystyle{{}##}$\hfil
   \tabskip=\centering&\kern-\displaywidth
        \rlap{\hbox{(\rm##\unskip)}}\tabskip=\displaywidth\crcr
               #1\crcr}}
\def\overrightarrow{\mathpalette\overrightarrow@}
\def\overrightarrow@#1#2{\vbox{\ialign{$##$\cr
    #1{-}\mkern-6mu\cleaders\hbox{$#1\mkern-2mu{-}\mkern-2mu$}\hfill
     \mkern-6mu{\to}\cr
     \noalign{\kern -1\p@\nointerlineskip}
     \hfil#1#2\hfil\cr}}}
\def\overleftarrow{\mathpalette\overleftarrow@}
\def\overleftarrow@#1#2{\vbox{\ialign{$##$\cr
     #1{\leftarrow}\mkern-6mu\cleaders\hbox{$#1\mkern-2mu{-}\mkern-2mu$}\hfill
      \mkern-6mu{-}\cr
     \noalign{\kern -1\p@\nointerlineskip}
     \hfil#1#2\hfil\cr}}}
\def\overleftrightarrow{\mathpalette\overleftrightarrow@}
\def\overleftrightarrow@#1#2{\vbox{\ialign{$##$\cr
     #1{\leftarrow}\mkern-6mu\cleaders\hbox{$#1\mkern-2mu{-}\mkern-2mu$}\hfill
       \mkern-6mu{\to}\cr
    \noalign{\kern -1\p@\nointerlineskip}
      \hfil#1#2\hfil\cr}}}
\def\underrightarrow{\mathpalette\underrightarrow@}
\def\underrightarrow@#1#2{\vtop{\ialign{$##$\cr
    \hfil#1#2\hfil\cr
     \noalign{\kern -1\p@\nointerlineskip}
    #1{-}\mkern-6mu\cleaders\hbox{$#1\mkern-2mu{-}\mkern-2mu$}\hfill
     \mkern-6mu{\to}\cr}}}
\def\underleftarrow{\mathpalette\underleftarrow@}
\def\underleftarrow@#1#2{\vtop{\ialign{$##$\cr
     \hfil#1#2\hfil\cr
     \noalign{\kern -1\p@\nointerlineskip}
     #1{\leftarrow}\mkern-6mu\cleaders\hbox{$#1\mkern-2mu{-}\mkern-2mu$}\hfill
      \mkern-6mu{-}\cr}}}
\def\underleftrightarrow{\mathpalette\underleftrightarrow@}
\def\underleftrightarrow@#1#2{\vtop{\ialign{$##$\cr
      \hfil#1#2\hfil\cr
    \noalign{\kern -1\p@\nointerlineskip}
     #1{\leftarrow}\mkern-6mu\cleaders\hbox{$#1\mkern-2mu{-}\mkern-2mu$}\hfill
       \mkern-6mu{\to}\cr}}}
\def\sqrt#1{\radical"270370 {#1}}
\def\dots{\relax\ifmmode\let\next=\ldots\else\let\next=\tdots@\fi\next}
\def\tdots@{\unskip\ \tdots@@}
\def\tdots@@{\futurelet\next\tdots@@@}
\def\tdots@@@{$\mathinner{\ldotp\ldotp\ldotp}\,
   \ifx\next,$\else
   \ifx\next.\,$\else
   \ifx\next;\,$\else
   \ifx\next:\,$\else
   \ifx\next?\,$\else
   \ifx\next!\,$\else
   $ \fi\fi\fi\fi\fi\fi}
\def\text{\relax\ifmmode\let\next=\text@\else\let\next=\text@@\fi\next}
\def\text@@#1{\hbox{#1}}
\def\text@#1{\mathchoice
 {\hbox{\everymath{\displaystyle}\def\textfonti{\the\textfont1 }%
    \def\textfontii{\the\textfont2 }\textdef@@ T#1}}
 {\hbox{\everymath{\textstyle}\def\textfonti{\the\textfont1 }%
    \def\textfontii{\the\textfont2 }\textdef@@ T#1}}
 {\hbox{\everymath{\scriptstyle}\def\textfonti{\the\scriptfont1 }%
   \def\textfontii{\the\scriptfont2 }\textdef@@ S\rm#1}}
 {\hbox{\everymath{\scriptscriptstyle}\def\textfonti{\the\scriptscriptfont1 }%
   \def\textfontii{\the\scriptscriptfont2 }\textdef@@ s\rm#1}}}
\def\textdef@@#1{\textdef@#1\rm \textdef@#1\bf
   \textdef@#1\sl \textdef@#1\it}

\def\textdef@#1#2{\def\next{\csname\expandafter\eat@\string#2fam\endcsname}%
\if S#1\edef#2{\the\scriptfont\next\relax}%
 \else\if s#1\edef#2{\the\scriptscriptfont\next\relax}%
 \else\edef#2{\the\textfont\next\relax}\fi\fi}
\scriptfont\itfam=\tenit \scriptscriptfont\itfam=\tenit
\scriptfont\slfam=\tensl \scriptscriptfont\slfam=\tensl
\mathcode`\0="0030
\mathcode`\1="0031
\mathcode`\2="0032
\mathcode`\3="0033
\mathcode`\4="0034
\mathcode`\5="0035
\mathcode`\6="0036
\mathcode`\7="0037
\mathcode`\8="0038
\mathcode`\9="0039
\def\Cal{\relax\ifmmode\let\next=\Cal@\else
     \def\next{\errmessage{Use \string\Cal\space only in math mode}}\fi\next}
\def\Cal@#1{{\fam2 #1}}
\def\bold{\relax\ifmmode\let\next=\bold@\else
   \def\next{\errmessage{Use \string\bold\space only in math
      mode}}\fi\next}\def\bold@#1{{\fam\bffam #1}}
\mathchardef\Gamma="0000
\mathchardef\Delta="0001
\mathchardef\Theta="0002
\mathchardef\Lambda="0003
\mathchardef\Xi="0004
\mathchardef\Pi="0005
\mathchardef\Sigma="0006
\mathchardef\Upsilon="0007
\mathchardef\Phi="0008
\mathchardef\Psi="0009
\mathchardef\Omega="000A
\mathchardef\varGamma="0100
\mathchardef\varDelta="0101
\mathchardef\varTheta="0102
\mathchardef\varLambda="0103
\mathchardef\varXi="0104
\mathchardef\varPi="0105
\mathchardef\varSigma="0106
\mathchardef\varUpsilon="0107
\mathchardef\varPhi="0108
\mathchardef\varPsi="0109
\mathchardef\varOmega="010A
\font\dummyft@=dummy
\fontdimen1 \dummyft@=\z@
\fontdimen2 \dummyft@=\z@
\fontdimen3 \dummyft@=\z@
\fontdimen4 \dummyft@=\z@
\fontdimen5 \dummyft@=\z@
\fontdimen6 \dummyft@=\z@
\fontdimen7 \dummyft@=\z@
\fontdimen8 \dummyft@=\z@
\fontdimen9 \dummyft@=\z@
\fontdimen10 \dummyft@=\z@
\fontdimen11 \dummyft@=\z@
\fontdimen12 \dummyft@=\z@
\fontdimen13 \dummyft@=\z@
\fontdimen14 \dummyft@=\z@
\fontdimen15 \dummyft@=\z@
\fontdimen16 \dummyft@=\z@
\fontdimen17 \dummyft@=\z@
\fontdimen18 \dummyft@=\z@
\fontdimen19 \dummyft@=\z@
\fontdimen20 \dummyft@=\z@
\fontdimen21 \dummyft@=\z@
\fontdimen22 \dummyft@=\z@
\def\fontlist@{\\{\tenrm}\\{\sevenrm}\\{\fiverm}\\{\teni}\\{\seveni}%
 \\{\fivei}\\{\tensy}\\{\sevensy}\\{\fivesy}\\{\tenex}\\{\tenbf}\\{\sevenbf}%
 \\{\fivebf}\\{\tensl}\\{\tenit}\\{\tensmc}}
\def\dodummy@{{\def\\##1{\global\let##1=\dummyft@}\fontlist@}}
\newif\ifsyntax@
\newcount\countxviii@
\def\newtoks@{\alloc@5\toks\toksdef\@cclvi}
\def\nopages@{\output={\setbox\z@=\box\@cclv \deadcycles=\z@}\newtoks@\output}
\def\syntax{\syntax@true\dodummy@\countxviii@=\count18
\loop \ifnum\countxviii@ > \z@ \textfont\countxviii@=\dummyft@
   \scriptfont\countxviii@=\dummyft@ \scriptscriptfont\countxviii@=\dummyft@
     \advance\countxviii@ by-\@ne\repeat
\dummyft@\tracinglostchars=\z@
  \nopages@\frenchspacing\hbadness=\@M}
\def\magstep#1{\ifcase#1 1000\or
 1200\or 1440\or 1728\or 2074\or 2488\or 
 \errmessage{\string\magstep\space only works up to 5}\fi\relax}
{\lccode`\2=`\p \lccode`\3=`\t 
 \lowercase{\gdef\tru@#123{#1truept}}}

\def\scaletype#1{\mag=#1\relax
 \hsize=\expandafter\tru@\the\hsize
 \vsize=\expandafter\tru@\the\vsize
 \dimen\footins=\expandafter\tru@\the\dimen\footins}

\def\scalefont#1#2\andcallit#3{\edef\font@{\the\font}#1\font#3=
  \fontname\font\space scaled #2\relax\font@}
\def\Mag@#1#2{\ifdim#1<1pt\multiply#1 #2\relax\divide#1 1000 \else
  \ifdim#1<10pt\divide#1 10 \multiply#1 #2\relax\divide#1 100\else
  \divide#1 100 \multiply#1 #2\relax\divide#1 10 \fi\fi}
\def\scalelinespacing#1{\Mag@\baselineskip{#1}\Mag@\lineskip{#1}%
  \Mag@\lineskiplimit{#1}}
\def\wlog#1{\immediate\write-1{#1}}
\catcode`\@=\active

\nopagenumbers
\TagsOnRight
\scaletype{\magstep1}
\scalelinespacing{\magstep1}
\pagewidth{15truecm}
\pageheight{22truecm}
\headline={\rightline{\folio}}
\voffset=2\baselineskip

\title Transversal fluctuations for increasing subsequences on the plane
\endtitle
\bigskip
\bigskip
\author Kurt Johansson \endauthor
\vskip 10truecm
\noindent
Version: Revised June 30, 1999
\newline Address: \newline Department of Mathematics\newline
Royal Institute of Technology\newline
S-100 44 Stockholm \newline
Sweden\newline
e-mail kurtj\@math.kth.se\newline
fax. +46 8 723 17 88

\newpage

\subheading{Abstract} Consider a realization of a Poisson process in
$\Bbb R^2$ with intensity 1 and take a maximal up/right path from the
origin to $(N,N)$ consisting of line segments between the points,
where maximal means that it contains as many points as possible. The
number of points in such a path has fluctuations of order $N^{\chi}$,
where $\chi=1/3$, [BDJ]. Here we show that typical deviations of a
maximal path from the diagonal $x=y$ is of order $N^{\xi}$ with
$\xi=2/3$. This is consistent with the scaling identity $\chi=2\xi-1$
which is believed to hold in many random growth models.
\bigskip\bigskip
\subheading{Mathematics Subject Classification}
\newline
Primary: 60K35, 82C24
\newline
Secondary: 82B24, 82B44

\newpage

\heading 1. Introduction and results\endheading
\medskip
The fluctuations in many random growth models, for example in
first-passage percolation, are described by two exponents, $\chi$ and
$\xi$, see e.g. [KS] and [LNP]. The exponent $\chi$ describes the
longitudinal whereas $\xi$ describes the transversal fluctuations. In
first-passage percolation the length of a minimizing path from the
origin to $(N,N)$ has fluctuations of order $N^{\chi}$, and the
minimizing path has typical deviations from the diagonal $x=y$ of order
$N^{\xi}$. General heuristic arguments (see [KS]) suggest that the
scaling identity $\chi=2\xi-1$ is valid in any dimension, compare the
heuristic argument below. In two dimensions it is predicted that
$\chi=1/3$ and hence we should have $\xi=2/3$. Since $\xi>1/2$ one
says that the minimizing path is {\it superdiffusive}.

We will consider a related model where it is known that $\chi=1/3$ and
prove that in this model we actually have $\xi=2/3$. The model is a
Poissonized version of the problem of the longest increasing
subsequence in a random permutation introduced in [Ha],
see also [AD]. In this model one considers a Poisson process with
intensity 1 in $\Bbb R_+^2$ and looks at a maximal up/right path from the
origin to $(N,N)$ consisting of line segments between the  Poisson points,
where maximal means that it contains as many points as possible. The
length of a path is the number of Poisson points in the path, and the
length of a maximal path has fluctuations of order $N^{1/3}$, see
[BDJ]. In this paper we will prove that the typical deviations of the
maximal paths from $x=y$ are of order $N^{2/3}$.

The proof uses the line of argument, for first-passage percolation
models, initiated in [NP], to prove
$\chi'\ge 2\xi-1$ (where $\chi'$ is closely related to $\chi$), and
[LNP] to prove lower (superdiffusive) bounds on a suitably defined
$\xi$. A related argument was
used to analyze the corresponding problem for crossing Brownian motion
in a Poissonian potential in [W\"u], and the present paper follows
the arguments in [W\"u]. A heuristic argument goes as follows. The
length of a typical maximal path from the origin to $(x,y)$ is $\sim
2\sqrt{xy}$, see [AD]. Hence, a maximal path from the origin to $(N,N)$
that passes through $(N(t-\delta),N(t+\delta))$, $0<t<1$, $\delta$ small,
is shorter by the amount
$$
2\sqrt{N(t-\delta)N(t+\delta)}+
2\sqrt{N(1-t+\delta)N(1-t-\delta)}-2\sqrt{N^2}.
$$
This should be of the same order as the length fluctuations,
i.e. $O(N^{\chi})$, which gives $\delta^2=O(N^{\chi-1})$. Thus,
$N^\xi\sim N\delta\sim N^{\chi/2+1/2}$, that is $2\xi-1=\chi$ and
hence $\xi=2/3$ since $\chi=1/3$. The argument used below essentially
makes this rigorous.

We will now give the precise definitions. Let $\Bbb P$ denote the
Poissonian law with fixed intensity 1 on the space $\Omega$ of locally
finite, simple, pure point measures on $\Bbb R^2$;
$\omega=\sum_i\delta_{\zeta_i}\in\Omega$, $\zeta_i=(x_i,y_i)$ are the
points in $\omega$. Write $(x,y)\prec (x',y')$ if $x<x'$ {it and}
$y<y'$. Given $\omega$ and two points $w\prec w'$ in $\Bbb R^2$ an
{\it up/right path} $\pi$ from $w$ to $w'$ is a subsequence
$\{\zeta_{i_k}\} _{k=1}^M$ of points in $\omega$ such that
$$
w\prec\zeta_{i_1}\prec\dots\prec\zeta_{i_M}\prec w'.
$$
The length, $|\pi|$, of $\pi$ is $M$, the number of Poisson points in
the path. Let $\Pi(w,w';\omega)$ denote the set of all up/right paths
from $w$ to $w'$ in $\omega$. If $K$ is a convex subset of $\Bbb R^2$
we let $\Pi^K(w,w';\omega)$ denote all up/right paths $\pi$ from $w$
to $w'$ inside $K$, i.e. $\pi\subseteq K$ and $w,w'\in K$. Let
$$
d(w,w';\omega)=\max\{|\pi|\,;\,\pi\in\Pi(w,w';\omega)\},
$$
and
$$
d^K(w,w';\omega)=\max\{|\pi|\,;\,\pi\in\Pi^K(w,w';\omega)\}.
$$

Let $\ell_N(\sigma)$ denote the length of a longest increasing
subsequence in a random permutation $\sigma\in S_N$ (uniform
distribution). If $i_1<\dots <i_n$ and $\sigma(i_1)<\dots<\sigma(i_n)$
we have an increasing subsequence of length $n$ and $\ell_N(\sigma)$
is the length of the longest such sequence. We define the
Poissonized distribution function by
$$
\phi_n(\lambda)=e^{-\lambda}\sum_{N=0}^\infty
\frac{\lambda^N}{N!}P[\ell_N(\sigma)\le n],
$$
[$\ell_0(\sigma)\equiv 0$]. Let $a(w,w')$ denote the area of the
rectangle $[w,w']$ with corners at $w$ and $w'$.
Now,
$$
\Bbb P[d(w,w')\le n]=\sum_{N=0}^\infty\Bbb P[d(w,w')\le
n\,\bigl|\,\omega([w,w'])=N] \Bbb P[\omega([w,w'])=N],
$$
and, see [Ha] or [AD], $\Bbb P[d(w,w')\le
n\,\bigl|\,\omega([w,w'])=N]=P[\ell_N(\sigma)\le n]$. Hence
$$
\Bbb P[d(w,w')\le n]=\phi_n(a(w,w')).\tag 1.1
$$
By Lemma 7.1 in [BDJ] we have a very good control of the function
$\phi_n(\lambda)$. Let 
$$
t=2^{1/3}(n+1)^{-1/3}(n+1-2\sqrt{\lambda}).\tag 1.2
$$
Then for any fixed $t$ in $\Bbb R$,
$$
\lim_{\lambda\to\infty}\phi_n(\lambda)=F(t),\tag 1.3
$$
where $F(t)$ is the Tracy-Widom largest eigenvalue distribution for
GUE, see [TW] and [BDJ]. 
The distribution function $F(t)$ is given by
$$
F(t)=\exp\bigl(-\int_t^\infty(x-t)u(x)^2dx\bigr),
$$
where $u(x)$ is the solution of the Painlev\'e II equation
$$
u''(x)=2u(x)^3+xu(x),\quad\text{and $u(x)\sim\text{Ai\,}(x)$ as
$x\to\infty$},
$$
where $\text{Ai\,}(x)$ is the Airy function. From this formula and the
asymptotics of $u(x)$, see [BDJ], it follows that $0<F(0)<1$, which
will be used below.
Furthermore we have the following estimates.
There are positive constants $\delta,T_0, c_1, c_2$ so that
if $T_0\le t\le 2^{-2/3}(n+1)^{2/3}$, then
$$
|\log\phi_n(\lambda)|\le  c_1\exp(-c_2t^{3/2}),\tag 1.4
$$
and if $-\delta(n+1)^{2/3}\le t\le -T_0$, then
$$
\phi_n(\lambda)\le c_1\exp(c_2t^3),\tag 1.5
$$
for all sufficently large $n$.
The estimate (1.4) also follows from the results in [Se]. These
estimates will be
important in the proof of our theorem.

Let $C(\gamma,N)$ be the cylinder of width $N^\gamma$ from 0 to
$w_N=(N,N)$: 
$$
C(\gamma,N)=\{(x,y)\,;\,0\le x+y\le 2N\,,\, -\sqrt{2}N^{\gamma}\le
-x+y\le \sqrt{2}N^{\gamma}\}.
$$
Denote by
$$
\Pi_{\max}(w,w';\omega)=\{\pi\in\Pi(w,w';\omega)\,;\,
|\pi|=d(w,w';\omega)\},
$$
the set of maximal paths from $w$ to $w'$. We are interested in the
size of the fluctuations of maximal paths around the diagonal $x=y$,
{\it the transversal fluctuations}. Let $A_N^\gamma$ be the event that all
maximal paths from $0$ to $w_N$ are contained in the cylinder
$C(\gamma,N)$,
$$
A_N^\gamma=\{\omega\in\Omega\,;\,\text{for all
$\pi\in\Pi_{\max}(0,w_N;\omega)$ we have $\pi\subseteq C(\gamma,N)$}\}.
$$
The {\it exponent of transversal fluctuations}, $\xi$, is then defined
by
$$
\xi=\inf\{\gamma >0\,;\,\liminf_{N\to\infty}\Bbb
P[A_N^\gamma]=1\}.\tag 1.6
$$
We can now state the main result of the paper.
\proclaim{Theorem 1.1} For the model defined above the exponent of
transversal fluctuations $\xi=2/3$.
\endproclaim

The proof of the theorem occupies the next section.

\proclaim{Remark 1.2}\rm We can consider the analogous problem for the
growth model introduced in [Jo]. Let $w(i,j)$, $(i,j)\in\Bbb Z_+^2$,
be independent geometrically (or exponentially) distributed random
variables and consider
$$
G(N)=\max\{\sum_{(i,j)\in\pi}w(i,j)\,;\,\text{$\pi$ an up/right path from
$(1,1)$ to $(N,N)$}\}.
$$
In [Jo] it is proved that there are positive constants $a$ and $b$ so
that $(G(N)-aN)/bN^{1/3}$ converges in distribution to a random
variable with distribution function $F(t)$. In analogy with above we can
consider the transversal deviations of a maximal path and define the
exponent $\xi$. If we had large deviation estimates for $\Bbb
P[G(N)\le n]$ analogous to (1.4) and (1.5) we could copy the proof
given in the next section and show that $\xi=2/3$ in this case
also. In [Jo] an estimate like (1.4) is proved, but (1.5) is
open. It follows from [BR] that 
$\Bbb P[G(N)\le n]$ is given by a certain $n\times n$ Toeplitz
determinant just as $\phi_n(\lambda)$, and it might be possible to
prove the analogue of (1.5) using Riemann-Hilbert techniques as in [BDJ].

\bigskip
\heading 2. Proof of $\xi\ge 2/3$\endheading
\medskip
We will first prove that $\xi\ge 2/3$. Pick $\gamma\in (\xi,1)$ and $\epsilon
>0$ (small). 
That $\xi<1$ follows from the proof in sect. 3 that $\xi\le 2/3$,
which is independent of the present section.
By the definition of $\xi$ there is an $N_0$ such that
$$
\Bbb P[A_N^\gamma]\ge 1-\epsilon\tag 2.1
$$
for all $N\ge N_0$. If $\omega\in A_N^\gamma$, then every maximal path from
$0$ to $w_N$ is contained in the cylinder $C(\gamma,N)$, so writing
$C_1=C(\gamma,N)$, we see that
$d^{C_1}(0,w_N;\omega)=d(0,w_N;\omega)$. Hence, by (2.1),
$$
\Bbb P[d^{C_1}(0,w_N)=d(0,w_N)]\ge 1-\epsilon,\tag 2.2
$$
if $N\ge N_0$.

Set
$\bold{v}_1=(1/\sqrt{2}, 1/\sqrt{2})$ and $\bold{v}_2=(-1/\sqrt{2},
1/\sqrt{2})$.
Let $m_N=3N^\gamma \bold{v}_2$ and let $C_2$ be the
cylinder $C_2=C_1+m_N$. Pick a $b$ such that $\gamma<b<1$, and assume
that $N$ is so large that $N^b-4N^\gamma>0$.
Define the points $A,B,C$ on the sides of $C_2$ by
$$\align
\overline{OA}&=(N^b+2N^\gamma)\bold{v}_1+2N^\gamma\bold{v}_2,\\
\overline{OB}&=(N^b+4N^\gamma)\bold{v}_1+4N^\gamma\bold{v}_2,\\
\overline{OC}&=N^b\bold{v}_1+4N^\gamma\bold{v}_2.
\endalign$$
$ABC$ is a right angle triangle with the right angle at $A$, the side
$AB$ is vertical with $A$ on the lower side of $C_2$ and $B$ on the
upper side. Divide the vertical side $AB$ into $K=K(N)$ segments
$z_{i-1}z_i$, $i=1,\dots ,K$, where $z_0=A$ and $z_K=B$. Let $L_i$ be
the part of the straight line through $z_i$, parallel to the
$x$-axis, lying in $C_2$. The parallelogram between $L_{i-1}$ and
$L_i$ in $C_2$ is denoted by $F_i$, $i=1,\dots,K$. We also define the
analogous geometrical objects at the other end of the cylinder, close
to $m_N+w_N$, by translating the whole picure by
$t_N=\sqrt{2}N-6N^\gamma-2N^b$, $z_i'=z_i+t_N\bold{v}_1$,
$F_i'=F_i+t_N\bold{v}_1$, $\overline{OA}'=\overline{OA}+t_N\bold{v}_1$ and
$\overline{OB}'=\overline{OB}+t_N\bold{v}_1$.

Given a Borel set $F$, $\omega(F)$ is the number of Poisson points in
$F$. Let $\pi=\{\zeta_1,\dots,\zeta_M\}$,
$\zeta_1\prec\dots\prec\zeta_M$, be a maximal path in
$\Pi^{C_2}(m_N,m_N+w_N;\omega)$ and let $\pi^\ast$ be the curve obtained
by joining $\zeta_i$ and $\zeta_{i+1}$, $i=0,\dots,K$, by straight
line segments, $\zeta_0=m_N$ and $\zeta_{K+1}=m_N+w_N$. The curve
$\pi^\ast$ intersects $AB$ at some point $P$ and $A'B'$ at some point
$Q$. The point $P$ belongs to $\bar F_i$ and $Q$ to $\bar F_j'$ for
some $i,j$. We will write $z(\omega)=z_i$ and $z'(\omega)=z_j'$. (If
$P=z_i$ for some $i$ we let $z(\omega)=z_i$ and analogously for $Q$.)
If we set $D_N(\omega)=\max_i\omega(\bar F_i)+\max_j\omega(\bar
F_j')$, then
$$\align
d^{C_2}(m_N,m_N+w_N)&\le
d^{C_2}(m_N,z(\omega))+d^{C_2}(z(\omega),z'(\omega))\\&+
d^{C_2}(z'(\omega),m_N+w_N)+D_N(\omega).\tag 2.3
\endalign$$
Note that $z(\omega)\in\Cal A\doteq\{z_0,\dots,z_K\}$ and
$z'(\omega)\in\Cal A'\doteq\{z'_0,\dots,z'_K\}$.
\proclaim{Lemma 2.1} Let $K=[8N^{2\gamma}]+1$. Then
$$
\Bbb P[D_N(\omega)\ge d]\le C(8N^{2\gamma}+1)e^{-d/2},\tag 2.4
$$
for all $d\ge 1$, where $C$ is a numerical constant.
\endproclaim

\demo{Proof} Since
$$
\{D_N(\omega)\ge d\}\subseteq\{\max_i\omega(\bar F_i)\ge\frac d2\}
\cup
\{\max_j\omega(\bar F'_j)\ge\frac d2\}
$$
we have
$$
\Bbb P[D_N(\omega )\ge d]\le 2K\Bbb P[\omega(\bar F_1)\ge d/2].\tag
2.5
$$
Here we have used the fact that all the random variables $\omega(\bar
F_i)$, $\omega(\bar F'_j)$ are identically distributed. The area of
$\bar F_1$ is $8N^{2\gamma}/K=\lambda$, and thus
$$
\Bbb P[\omega(\bar F_1)\ge d/2]\le\sum_{j=[d/2]}^\infty
e^{-\lambda}\frac{\lambda^j}{j!}\le C\sum_{j=[d/2]}^\infty e^{-\lambda
f(j/\lambda)} ,\tag 2.6
$$
where $C$ is a numerical constant and $f(x)=x\log x+1-x$. Here we have
used Stirling's formula. Note that $f(x)\ge x$ if $x\ge 9$ say. Choose
$K=[8N^{2\gamma}]+1$, so that $\lambda\le 1$, and assume that $d\ge
18$. Then, by (2.6),
$$
\Bbb P[\omega(\bar F_1)\ge d/2]\le C\sum_{j=[d/2]}^\infty e^{-j}\le
Ce^{-d/2}
$$
and introducing this estimate into (2.5) yields
$$
\Bbb P[\omega(\bar F_1)\ge d]\le C(1+8N^{2\gamma})e^{-d/2}
$$
for all $N\ge 1$, $d\ge 1$.
\enddemo
\rightline{Q.E.D}
It follows from the estimate (2.4) that
$$
\Bbb P[D_N(\omega)\le 5\log N]\ge 1-\epsilon,\tag 2.7
$$
for all sufficiently large $N$. 

Next, choose $\kappa_1$ and $\kappa_2$
so that $0<\kappa_1<1/3<\kappa_2<1$.

\proclaim{Lemma 2.2} Assume that (2.1) holds. There is a numerical
constant $\eta\in (0,1)$, such that if $\epsilon\le\eta$ and $N$ is
sufficiently large, then
$$
\Bbb P[d^{C_1}(0,w_N)-d^{C_2}(m_N,m_N+w_N)\le
-N^{\kappa_1}]\ge\eta.\tag 2.8
$$
Furthermore, for $N$ sufficiently large,
$$
\Bbb P[|d(0,z(\omega))-2\sqrt{a(0,z(\omega))}|\le N^{b\kappa_2}]\ge
1-\epsilon, \tag 2.9
$$
$$
\Bbb P[|d(m_N,z(\omega))-2\sqrt{a(m_N,z(\omega))}|\le N^{b\kappa_2}]\ge
1-\epsilon, \tag 2.10
$$
$$
\Bbb P[|d(z'(\omega),w_N)-2\sqrt{a(z'(\omega),w_N)}|\le N^{b\kappa_2}]\ge
1-\epsilon, \tag 2.11
$$
$$
\Bbb P[|d(z'(\omega),w_N+m_N)-2\sqrt{a(z'(\omega),w_N+m_N)}|
\le N^{b\kappa_2}]\ge
1-\epsilon, \tag 2.12
$$
\endproclaim
\demo{Proof} The random variables $d^{C_1}(0,w_N)$ and
$d^{C_2}(m_N,m_N+w_N)$ are independent. Thus
$$\align
&\Bbb P[d^{C_1}(0,w_N)-d^{C_2}(m_N,m_N+w_N)\le
-N^{\kappa_1}]\\
&\ge\Bbb P[d^{C_1}(0,w_N)-2N\le
0\,\,\text{and}\,\,
d^{C_2}(m_N,m_N+w_N)-2N\ge N^{\kappa_1}]\\
&=\Bbb P[d^{C_1}(0,w_N)-2N\le 0]\cdot
\Bbb P[d^{C_1}(0,w_N)-2N\ge N^{\kappa_1}].\tag 2.13
\endalign
$$
If $\omega\in A_N^\gamma$, then $d^{C_1}(0,w_N)=d(0,w_N)$, and
consequently the last expression in (2.13) is greater than or equal to
$$\align
&\Bbb P[\{d(0,w_N)-2N\le 0\}\cap A_N^\gamma]\cdot
\Bbb P[\{d(0,w_N)-2N\ge N^{\kappa_1}\}\cap A_N^\gamma]
\\
&\ge
(\Bbb P[d(0,w_N)-2N\le 0]+\Bbb P[A_N^\gamma]-1)\\&\times
(\Bbb P[d(0,w_N)-2N\ge N^{\kappa_1}]+\Bbb
P[A_N^\gamma]-1).
\tag 2.14
\endalign
$$
By (1.1),
$$
\Bbb P[d(0,w_N)-2N\le 0]=\phi_{2N}(N^2).
$$
It follows from (1.3) that $\phi_{2N}(N^2)\to F(0)$ as
$N\to\infty$. Furthermore, since $\kappa_1<1/3$,
$$
\Bbb P[d(0,w_N)-2N\ge N^{\kappa_1}]=1-\phi_{[2N+N^{\kappa_1}]}(N^2)\to 1-F(0),
$$
as $N\to\infty$, again by (1.3) and the fact that $\phi_n(\lambda)$ is
increasing in $n$. 
Let $\eta=\frac 13F(0)(1-F(0)>0$. If $N$ is sufficiently large then
$\Bbb P[d(0,w_N)-2N\le 0]\ge F(0)-\eta$ and 
$\Bbb P[d(0,w_N)-2N\ge N^{\kappa_1}]\ge 1-F(0)-\eta$. Since $\Bbb
P[A_N^\gamma]\ge 1-\epsilon$ by (2.1) we see that the right hand side
of (2.14) is $\ge (F(0)-2\eta)(1-F(0)-2\eta)\ge \eta$. 
This proves (2.8).

Next, we will prove (2.9). The proofs of (2.10), (2.11) and (2.12) are
completely analogous. Note that
$$
\Bbb P[|d(0,z(\omega))-2\sqrt{a(0,z(\omega))}|\le N^{b\kappa_2}]\ge
\Bbb P[\bigcap_{j=0}^K\{|d(0,z_j)-
2\sqrt{a(0,z_j)}|\le N^{b\kappa_2}\}]
$$
so it suffices to show that
$$
\sum_{j=1}^K\Bbb P[|d(0,z_j)-
2\sqrt{a(0,z_j)}|\ge N^{b\kappa_2}]\le\epsilon\tag 2.15
$$
for all sufficiently large $N$. We have $z_j=\frac
1{\sqrt{2}}(N^b,N^b+r_j)$, where $4N^\gamma\le r_j\le 8N^\gamma$, so
$a(0,z_j)=\frac 12(N^{2b}+N^br_j)\doteq a_j$. Now,
$$
\Bbb P[d(0,z_j)-2\sqrt{a_j}\le
-N^{b\kappa_2}]=\phi_{[2\sqrt{a_j}-N^{b\kappa_2}]}(a_j).
$$
In this case $t$ defined by (1.2) is $\sim-2^{1/6}N^{b\kappa_2-b/3}$ and
since $1/3<\kappa_2<1$, the condition for (1.5) is fulfilled if $N$ is
sufficiently large and we get
$$
\Bbb P[d(0,z_j)-2\sqrt{a_j}\le -N^{b\kappa_2}]\le
c_3\exp (-c_4N^{3b(\kappa_2-1/3)})\tag 2.16
$$
for some positive constants $c_3, c_4$ and all $j$. Similarly we can
use (1.4) to prove that 
$$
\Bbb P[d(0,z_j)-2\sqrt{a_j}\ge N^{b\kappa_2}]\le c_5\exp (-c_6N^{\frac
32 b(\kappa_2-1/3)})\tag 2.17
$$
for some positive constants $c_5, c_6$ if $N$ is sufficiently
large. Using (2.16) and (2.17) we see that (2.15) holds if $N$
is sufficiently large since $K=[8N^{2\gamma}]+1$. This completes the
proof of the lemma.
\enddemo
\rightline{Q.E.D}
Denote by $B_N^\gamma$ the set of $\omega$ that satisfy all the
inequalities inside $\Bbb P[\phantom{X}]$ in (2.7) - (2.12). 
Then, by (2.7) and
Lemma 2.2,
$$
\Bbb P[B_N^\gamma]\ge\eta-5\epsilon.\tag 2.18
$$
Note that for any $\omega$,
$$
d(0,w_N)\ge
d(0,z(\omega))+d(z(\omega),z'(\omega))+d(z'(\omega),w_N).\tag 2.19
$$
The inequalities (2.3) and (2.19) give
$$
\align
&d^{C_2}(m_N,m_N+w_N)-d(0,w_N)
\le d^{C_2}(m_N,z(\omega))+d^{C_2}(z'(\omega),m_N+w_N)\\&-d(0,z(\omega))-
d(z'(\omega),w_N)+D_N(\omega).\tag 2.20
\endalign$$
Now, using (2.20), we see that for $\omega\in B_N^\gamma$, 
$$
\align
&d^{C_1}(0,w_N)-d(0,w_N)\\
&=
d^{C_1}(0,w_N)-d^{C_2}(m_N,m_N+w_N)+d^{C_2}(m_N,m_N+w_N)-d(0,w_N)\\
&\le -N^{\kappa_1}+4N^{b\kappa_2}+2\sqrt{a(m_N,z(\omega))}+
2\sqrt{a(z'(\omega),m_N+w_N)}\\&-2\sqrt{a(0,z(\omega))}-
2\sqrt{a(z'(\omega),w_N)}+5\log N.\tag 2.21
\endalign
$$
To proceed we need the following purely geometric lemma.

\proclaim{Lemma 2.3} For all sufficiently large $N$,
$$
\sqrt{a(m_N,z)}-\sqrt{a(0,z)}\le 10N^{2\gamma-b}\tag 2.22
$$
for any $z\in\Cal A$ and
$$
\sqrt{a(z',w_N+m_N)}-\sqrt{a(z',w_N)}\le 10N^{2\gamma-b}\tag 2.23
$$
for any $z'\in\Cal A'$.
\endproclaim
\demo{Proof} We will prove (2.22). The inequality (2.23) then follows
by symmetry. Now, $a(m_N,z_j)=(N^{b}+3N^\gamma)(N^b-3N^\gamma+r_j)/2$, 
$a(0,z_j)=(N^{2b}+r_jN^b)/2$ and hence
$$
\sqrt{a(m_N,z)}-\sqrt{a(0,z)}=\frac{a(m_N,z)-a(0,z)}
{\sqrt{a(m_N,z)}+\sqrt{a(0,z)}}
\le\frac{3r_jN^\gamma}{2\sqrt{2}N^b}\le 10N^{2\gamma-b}, 
$$
since $r_j\le 8N^\gamma$.
\enddemo
\rightline{Q.E.D.}

Introducing the estimates (2.22) and (2.23) into (2.21) we obtain
$$
d^{C_1}(0,w_N)-d(0,w_N)\le
-N^{\kappa_1}+5N^{b\kappa_2}+40N^{2\gamma-b}
$$
for all $\omega\in B_N^\gamma$ if $N$ is sufficiently large. Thus, by (2.18),
$$
\Bbb P[d^{C_1}(0,w_N)-d(0,w_N)\le
-N^{\kappa_1}+5N^{b\kappa_2}+40N^{2\gamma-b} ]
\ge\eta-5\epsilon\ge\frac{\eta}2,\tag 2.24
$$
if $\epsilon<\eta/10$ and $N$ is sufficiently large. But we also
have the estimate (2.2). These estimates are consistent for large $N$
only if
$$
\kappa_1\le\max\{b\kappa_2,2\gamma-b\}.\tag 2.25
$$
In this inequality we can let $\kappa_1\nearrow 1/3$ and
$\kappa_2\searrow 1/3$ to get $1/3\le\max\{b/3,2\gamma-b\}$ and since
$b<1$, we must have $1/3\le 2\gamma-b$. Here we can let
$\gamma\searrow\xi$ and $b\nearrow 1$ to get $1/3\le 2\xi-1$,
i.e. $\xi\ge 2/3$.
\bigskip
\heading 3. Proof of $\xi\le 2/3$\endheading
We turn now to the proof of the opposite inequality $\xi\le 2/3$. By
the definition (1.6) of $\xi$ we see that we have to show that if
$\gamma>2/3$, then
$$
\lim_{N\to\infty}\Bbb P[\Omega\setminus A_N^\gamma]=0.\tag 3.1
$$
If $\omega\in \Omega\setminus A_N^\gamma$, then there is a path
$\pi_0\in \Pi_{\max}(0,w_N;\omega)$ such that $\pi_0$ is not contained
in $C(\gamma,N)$. We take one such path. Fix $\gamma\in (2/3,1)$. Let
$\pi_0^\ast$ be the curve associated to $\pi_0$. Then $\pi_0^\ast$
intersects the upper and/or the lower sides of $C(\gamma,N)$. Assume
that it intersects the upper side. Define a sequence of points on the
upper side of $C(\gamma,N)$, $z_j=(jM/K,jM/K+\sqrt{2}N^\gamma)$, $0\le
j\le K$, where $M=N-\sqrt{2}N^\gamma$ and
$K=[2\sqrt{2}N^{1+\gamma}]+1$. Let $D_j$ be the parallelogram with
corners at $z_{j-1}$, $z_j$, $(jM/K,jM/K-\sqrt{2} N^\gamma)$ and 
$((j-1)M/K, (j-1)M/K-\sqrt{2}N^\gamma)$, $1\le j\le K$.

The curve $\pi_0^\ast$ intersects the upper side for the first time,
going from $0$ to $w_N$, in the line segment $z_{j-1}z_j$ for some
$j$. We set $z(\omega)=z_{j-1}$. By the choice of $z(\omega)$ we have
that
$$
d(0,w_N)\le d(0,z(\omega))+d(z(\omega),w_N)+\max_{1\le j\le
K}\omega(D_j).\tag 3.2
$$

In the case when $\pi_0^\ast$ does not intersect the upper side but
only the lower side, there is a last time where it intersects the
lower side and we can assign a point $z(\omega)$ on the lower side so
that (3.2) holds. This case is the image under the map $T_N:(x,y)\to
(N-x,N-y)$ of the first case. Let $\Cal C=\{z_j\}_{j=0}^K$ and let
$\Cal C'$ be the image of $\Cal C$ under $T_N$. 

\proclaim{Lemma 3.1} Set 
$$
\Lambda_N=\{\omega\,;\,\max_{1\le j\le K}\omega(D_j)\le 2\log N\},
$$
and for each $z\in\Cal C\cup\Cal C'$, $\delta\in(1/3,2\gamma-1)$,
$$\align
E_z=\{\omega\,&;\,d(0,z)\le 2\sqrt{a(0,z)}+a(0,z)^{\delta/2}+N^\delta\,\,
\\&\text{and}\,\, d(z,w_N)\le
2\sqrt{a(z,w_N)}+a(z,w_N)^{\delta/2}+N^\delta\}.
\endalign$$
For any given $\epsilon >0$, there is an $N_0$ such that if
$N\ge N_0$, then
$$
\Bbb P\biggl[\bigcup_{z\in\Cal C\cup\Cal C'}(\Omega\setminus
E_z)\cup(\Omega\setminus\Lambda_N)\biggr] \le\epsilon.\tag 3.3
$$
\endproclaim
\demo{Proof} An argument analogous to the one used in the proof of
Lemma 2.1 shows that there is a numerical constant $C$ so that
$$
\Bbb P[\Omega\setminus \Lambda_N]\le CN^{\gamma-1}.
$$
We consider $z\in\Cal C$, the case $z\in\Cal C'$ is analogous by
symmetry. Recall that $[z,w]$ denotes the rectangle with corners at $z$ and
$w$. If $a(0,z)\le N^{\delta/2}$, then 
$\Bbb P[\omega([0,z])\ge N^{\delta}]\le C\exp(-N^\delta/2)$ for
some numerical constant $C$, by Chebyshev's inequality. Since we trivially have
$d(0,z;\omega)\le\omega([0,z])$, we obtain
$$
\Bbb P[d(0,z)>2\sqrt{a(0,z)}+a(0,z)^{\delta/2}+N^\delta]\le
C\exp(-N^{\delta}/2) ,\tag 3.4
$$
provided $a(0,z)\le N^{\delta/2}$. Now, with $a=a(0,z)$,
$$
\Bbb P[d(0,z)>2\sqrt{a}+a^{\delta/2}+N^\delta]\le 1-\phi_{[2\sqrt{a}+
a^{\delta/2}]}(a).
$$
This last expression can be estimated using (1.4), which gives
$$
1-\phi_{[2\sqrt{a}+ a^{\delta/2}]}(a)\le
c_1'\exp(-c_2'a^{(\delta-1/3)/2}).
$$
If $a\ge N^{\delta/2}$, the right hand side is
$\le c_1'\exp(-c_2'N^{\delta (\delta -1/3)/4})$ and thus
$$
\Bbb P[d(0,z)>2\sqrt{a}+a^{\delta/2}+N^\delta]\le 
c_1'\exp(-c_2'N^{\delta (\delta -1/3)/4}).\tag 3.5
$$
We can prove estimates analogous to (3.4) and (3.5) with
$d(0,z)$ replaced by $d(z,w_N)$ in the same way. Bringing everything
together we see that (3.3) holds if $N$ is sufficiently large.
The lemma is proved.
\enddemo
\rightline{Q.E.D.}

Set 
$$
B_N^\gamma=(\Omega\setminus A_N^\gamma)\cap(\bigcap_{z\in\Cal C\cup\Cal C'} E_z
)\cap\Lambda_N.
$$
By Lemma 3.1, for $N\ge N_0$,
$$
\Bbb P[\Omega\setminus A_N^\gamma]\le \epsilon+\Bbb P[B_N^\gamma].\tag
3.6
$$
Since $a(0,z)\le N^2$ and $a(z,w_N)\le N^2$ for any $z\in\Cal
C\cup\Cal C'$, we see from (3.2) that for $\omega\in B_N^\gamma$,
$$
d(0,w_N)\le 2\log
N+4N^\delta+2\sqrt{a(0,z(\omega))}+\sqrt{a(z(\omega),w_N)}.\tag 3.7
$$
We need one more geometric lemma.

\proclaim{Lemma 3.2} For any $z\in\Cal C\cup\Cal C'$,
$$
\sqrt{a(0,z)}+\sqrt{a(z,w_N)}-\sqrt{a(0,w_N)}\le -N^{2\gamma-1},
\tag 3.8
$$
if $N$ is sufficiently large.
\endproclaim
\demo{Proof} Again, by symmetry, it suffices to consider the case
$z\in\Cal C$. Now,
$a(0,z_j)=j\frac MK (j\frac MK +\sqrt{2}N^\gamma)$ and $a(z_j,w_N)= 
(N-j\frac MK )(N-j\frac MK -\sqrt{2}N^\gamma)$. where $1\le j\le
K=[2\sqrt{2}N^{1+\gamma}] +1$ and $M=N-\sqrt{2}N^{\gamma}$.
Write $x=jM/KN$ and  $y=\sqrt{2} N^{\gamma-1}$, so that
$0\le x\le 1-y$. Then,
$$
\sqrt{a(0,z)}+\sqrt{a(z,w_N)}-\sqrt{a(0,w_N)}=Nf(x,y),\tag 3.9
$$
where
$$
f(x,y)=\sqrt{x^2+xy}+\sqrt{(1-x)^2-(1-x)y}.
$$
For a fixed $y\in (0,1)$ this function assumes its maximum in
$[0,1-y]$ at $x=(1-y)/2$, which gives $f(x,y)\le -y^2/2$. Inserting
this estimate into (3.9) and taking $y=\sqrt{2}N^{\gamma-1}<1$, which
is true if $N$ is large enough,
proves the lemma.
\enddemo
\rightline{Q.E.D.}

Combining the estimates (3.7) and (3.8), we see that 
$$
\Bbb P[B_N^\gamma]\le\Bbb P[d(0,w_N)-2\sqrt{a(0,w_N)}\le 2\log
N+4N^\delta-2N^{2\gamma-1}]. \tag 3.10
$$
To finish the proof we need
\proclaim{Lemma 3.3} If $\delta\in (1/3,2\gamma -1)$, $\gamma >2/3$,
then
$$
\lim_{N\to\infty}\Bbb P[d(0,w_N)-2\sqrt{a(0,w_N)}\le 2\log
N+4N^\delta-2N^{2\gamma-1}]=0.\tag 3.11
$$
\endproclaim
\demo{Proof} Since $\delta<2\gamma -1$, we have that $2\log
N+4N^\delta-2N^{2\gamma -1}\le -N^{2\gamma -1}$ if $N$ is
sufficiently large. Thus, by (1.1),
$$\align
\Bbb P[d(0,w_N)\le 2N+2\log
N+4N^\delta-2N^{2\gamma-1}]&\le 
\Bbb P[d(0,w_N)\le 2N
-N^{2\gamma-1}]\\&=\phi_{[2N-N^{2\gamma -1}]}(N^2).
\endalign$$
The identity (1.2) with $n=[2N-N^{2\gamma -1}]$ and $\lambda=N^2$
gives $t\sim -N^{2\gamma-4/3}$, and hence (1.5) gives us the
estimate
$$
\phi_{[2N-N^{2\gamma -1}]}(N^2)\le c_1\exp(-c_2'N^{6\gamma-4}),
$$
where $c_2'>0$. This proves the lemma.
\enddemo
\rightline{Q.E.D.}

Combining (3.11) with (3.6) and (3.10) we have proved (3.1). Thus
$\xi\le 2/3$ and we are done.

\medskip
\heading Acknowledgement\endheading
I thank M. W\"uthrich for explaining his work on the exponents $\chi$
and $\xi$. I also thank A. -S. Sznitman for inviting me to the
Forschungsinsitut f\"ur Mathematik at ETH, Z\"urich, where this work
was begun, and E. Rains and J. Baik 
for keeping me informed about the work [BR]. 
This work was supported by the Swedish Natural Science
Research Council (NFR).

\bigskip
\define\kl#1{\medskip\item{[#1]    }}

\flushpar \heading REFERENCES\endheading \medskip

\kl{AD} D. Aldous and P. Diaconis, \it Hammersley's Interacting
Particle Process and Longest Increasing Subsequences\rm , Prob. Th. and Rel.
Fields, {\bf 103}, (1995), pp. 199 - 213

\kl{BDJ} J. Baik, P. A. Deift and K. Johansson, \it On the distribution
of the longest increasing subsequence in a random permutation\rm, 
math.CO/98101105, to appear in J. Amer. Math. Soc.

\kl{BR} J. Baik and E. Rains, \it Algebraic aspects of increasing
subsequences \rm, \newline math.CO/ 9905083

\kl{Ha} J. M. Hammersley, \it A few seedlings of research\rm ,  In
\it Proc. Sixth Berkeley Symp. Math. Statist. and Probability\rm , Volume 1,
pp. 345 - 394, University of California Press, 1972

\kl{Jo} K. Johansson, \it Shape fluctuations and random matrices\rm,
math.CO/9903134

\kl{KS} J. Krug, H. Spohn, \it Kinetic Roughening of Growing
Interfaces\rm, in Solids far from Equilibrium: Growth, Morphology and
Defects , Ed. C. Godr\`eche, 479 - 582, Cambridge University Press,
1992

\kl{LNP} C. Licea, C. M. Newman and M. S. T. Piza, \it
Superdiffusivity in first-passage percolation\rm, Probab.Theory
Relat. Fields, {\bf 106}, (1996), 977 - 1005

\kl{NP} C. M. Newman, M. S. T. Piza, \it Divergence of Shape
Fluctuations in Two Dimensions\rm, Ann. Prob., {\bf 23}, (1995), 977 -
1005 

\kl{Se} T. Sepp\"al\"ainen, \it Large Deviations for Increasing
Sequences on the Plane\rm, \linebreak Probab. Theory Relat. Fields, {\bf 112},
(1998), 221 - 244

\kl{TW} C. A. Tracy, H. Widom, \it Level Spacing Distributions and
the Airy Kernel\rm, Commun. Math. Phys., {\bf 159}, (1994), 151 - 174

\kl{W\"u} M. V. W\"uthrich, \it Scaling identity for crossing Brownian
motion in a Poissonian potential\rm, Probab. Theory Relat. Fields,
{\bf 112}, (1998), 299 - 319

\end